\newcommand{\al}{\alpha}               
\newcommand{\be}{\beta}
\newcommand{\ga}{\gamma}               

\newcommand{\de}{\delta}               
\newcommand{\De}{\Delta}

\newcommand{\veps}{\varepsilon}        
\newcommand{\vphi}{\varphi}

\newcommand{\cal}{\mathcal}

\newcommand{\calf}{{\cal F}}

\newcommand{\calv}{{\cal V}}

\newcommand{\Dom}{{\rm Dom}}
\newcommand{\Fix}{{\rm Fix}}

\newcommand{\incl}{\subseteq}

\newcommand{\es}{\emptyset}          
\newcommand{\sm}{\setminus}
      
\newcommand{\limpl}{\Longrightarrow}

\newcommand{\oo}{\infty}

\def\dtends   {\stackrel {\it d}{\longrightarrow}}

\def\(V)tends {\stackrel {(\calv)}{\longrightarrow}}

\newcommand{\barr}{\begin{array}}        
\newcommand{\earr}{\end{array}}
\newcommand{\bcor}{\begin{corollary}}    
\newcommand{\ecor}{\end{corollary}}
\newcommand{\ben}{\begin{enumerate}}     
\newcommand{\een}{\end{enumerate}}
\newcommand{\beq}{\begin{equation}}       
\newcommand{\eeq}{\end{equation}}
\newcommand{\bex}{\begin{example}}        
\newcommand{\eex}{\end{example}}
\newcommand{\bit}{\begin{itemize}}        
\newcommand{\eit}{\end{itemize}}
\newcommand{\blemma}{\begin{lemma}}       
\newcommand{\elemma}{\end{lemma}}
\newcommand{\bproof}{\begin{proof}}       
\newcommand{\eproof}{\end{proof}}
\newcommand{\bprop}{\begin{proposition}}  
\newcommand{\eprop}{\end{proposition}}
\newcommand{\brem}{\begin{remark}}        
\newcommand{\erem}{\end{remark}}
\newcommand{\btab}{\begin{tabular}}       
\newcommand{\etab}{\end{tabular}}
\newcommand{\btheorem}{\begin{theorem}}   
\newcommand{\etheorem}{\end{theorem}}

\documentclass[reqno]{amsart}

\newtheorem{theorem}{\bf Theorem}
\newtheorem{corollary}{\bf Corollary}
\newtheorem{example}{\bf Example}
\newtheorem{lemma}{\bf Lemma}
\newtheorem{proposition}{\bf Proposition}
\newtheorem{remark}{\bf Remark}

\begin{document}

\title
[Wardowski Implicit Contractions in Metric Spaces]
{WARDOWSKI IMPLICIT CONTRACTIONS \\
IN METRIC SPACES}

\author{Mihai Turinici}
\address{
"A. Myller" Mathematical Seminar;
"A. I. Cuza" University;
700506 Ia\c{s}i, Romania
}
\email{mturi@uaic.ro}




\begin{abstract}
Most of the implicit contractions introduced by
Wardowski [Fixed Point Th. Appl., 2012, 2012:94]
are Matkowski type contractions.
\end{abstract}

\maketitle

\section{Introduction}
\setcounter{equation}{0}

Let $X$ be a nonempty set.
Call the subset $Y$ of $X$, 
{\it almost singleton} (in short: {\it asingleton})
provided $y_1,y_2\in Y$ implies $y_1=y_2$;
and {\it singleton},
if, in addition, $Y$ is nonempty;
note that, in this case,
$Y=\{y\}$, for some $y\in X$. 
Further, let $d:X\times X\to R_+:=[0,\oo[$
be a {\it metric} over $X$; 
the couple $(X,d)$ will be referred to as a
{\it metric space}.
Define a sequential $d$-convergence $(\dtends)$ 
on $X$, according to:
for each sequence $(x_n; n\ge 0)$ in $X$ 
and each $x\in X$,
$x_n \dtends x$ iff $d(x_n,x)\to 0$; i.e.:
\ben
\item[]
$\forall \veps> 0$, $\exists i=i(\veps)$:\
$i\le n$ $\limpl$ $d(x_n,x)\le \veps$;
\een
referred to as:
$x$ is the {\it d-limit} of $(x_n; n\ge 0)$.
Denote by $\lim_n(x_n)$ the set of all such elements;
if it is nonempty, then $(x_n; n\ge 0)$ is called
{\it $d$-convergent}.
The class of all $d$-convergent sequences will 
be indicated as Conv$(X,d)$;
note that, it is {\it separated}:
$\lim_n(x_n)$ is an asingleton, for each
sequence $(x_n)$ of $X$.
Further,
call the sequence $(x_n; \ge 0)$ in $X$, {\it $d$-Cauchy}, provided
$d(x_m,x_n)\to 0$ as $m,n\to \oo$, $m< n$; i.e.
\ben
\item[]
$\forall \veps> 0$, $\exists j=j(\veps)$:\
$j\le m< n$ $\limpl$ $d(x_m,x_n)\le \veps$.
\een
The class of all these will be indicated as Cauchy$(X,d)$;
it includes all sequences $(x_n; n\ge 0)$
that are {\it $d$-telescopic-Cauchy}
(in short: {\it $d$-tele-Cauchy}), introduced as
\ben \item[]
$\sum_n d(x_n,x_{n+1})(=d(x_0,x_1)+d(x_1,x_2)+...)< \oo$.
\een
By definition, (Conv$(X,d)$, Cauchy$(X,d)$) will be
called the conv-Cauchy structure attached to $(X,d)$.
Note that (as $d$=metric), $(X,d)$ is {\it regular}: any 
$d$-convergent sequence is $d$-Cauchy;
if the reciprocal holds too,
then $(X,d)$ is called {\it complete}.

Having these precise, take some $T\in \calf(X)$.
[Here, given the nonempty sets $A$  and $B$, 
$\calf(A,B)$ stands for 
the class of all functions $f:A\to B$; 
when $A=B$, we write $\calf(A,A)$ as $\calf(A)$].
Denote $\Fix(T):=\{z\in X; z=Tz\}$; any point of it will be called
{\it fixed} under $T$.
These points are to be determined in the context below,
comparable with the one in 
Rus \cite[Ch 2, Sect 2.2]{rus-2001}:

{\bf 1a)} 
We say that $T$ is a 
{\it Picard (resp., tele-Picard) operator} (modulo $d$) if,
for each $x\in X$, $(T^nx; n\ge 0)$ is 
$d$-Cauchy (resp., $d$-tele-Cauchy) 
and
$d$-convergent

{\bf 1b)} 
We say that $T$ is a 
{\it strong Picard (resp., tele-Picard) operator} (modulo $d$) if,
it is a Picard (resp., tele-Picard) operator (modulo $d$), 
and $\lim_n(T^nx)$ belongs to $\Fix(T)$

{\bf 1c)}
We say that $T$ is a 
{\it globally strong Picard (resp., tele-Picard) operator} (modulo $d$) if
it is a strong Picard (resp., tele-Picard) operator (modulo $d$),
and
$\Fix(T)$ is an asingleton (hence, a singleton).

Sufficient conditions for such properties 
will be stated in the class of
"functional" metric contractions. 
Call $T$, {\it $(d;\vphi)$-contractive} 
(for some $\vphi\in \calf(R_+)$), when
\ben
\item[] (a01)\ \ 
$d(Tx,Ty)\le  \vphi(d(x,y))$,\ $\forall x,y\in X$, $x\ne y$.
\een
The functions to be considered here are as follows.
Let us say that  $\vphi\in \calf(R_+)$ is {\it increasing}, in case
[$t_1\le t_2$ implies $\vphi(t_1)\le \vphi(t_2)$];
the class of all these will be denoted as $\calf(in)(R_+)$.
Further, call  $\vphi\in \calf(in)(R_+)$, 
{\it regressive} in case
$\vphi(0)=0$ and [$\vphi(t)< t$, $\forall t> 0$];
the subclass of all these will be denoted as $\calf(in,re)(R_+)$.
Finally, we shall say that $\vphi\in \calf(in,re)(R_+)$ 
is {\it Matkowski admissible} provided 
\ben
\item[]
$\vphi^n(t)\to 0$ as $n\to \oo$, for all $t\in R_+$;
\een
and {\it Matkowski tele-admissible}, in case
\ben
\item[]
$\Phi(t):=\sum_n \vphi^n(t)\
(=t+\vphi(t)+\vphi^2(t)+...) < \oo$,\ for all $t\in R_+$.
\een
[Here, for each $n\ge 0$, 
$\vphi^n$ stands for the $n$-th iterate of $\vphi$].
Clearly, any Matkowski tele-admissible function is Matkowski admissible;
the reciprocal is not in general true.
The following fixed point result in Matkowski \cite{matkowski-1975}
is our starting point.

\btheorem \label{t1}
Supposed that $T$ is $(d;\vphi)$-contractive,
for some increasing regressive function $\vphi\in \calf(in,re)(R_+)$.
In addition, let $(X,d)$ be complete. Then, 

{\bf i)}
If (in addition), $\vphi$ is Matkowski admissible, then 
$T$ is globally strong Picard (modulo $d$); precisely,
\beq \label{101}
\mbox{
$\Fix(T)=\{z\}$ and $T^nx\dtends z$, for each $x\in X$
}
\eeq

{\bf ii)}
If, moreover,
$\vphi$ is Matkowski tele-admissible,
then $T$ is globally strong tele-Picard (modulo $d$);
hence, the above relation holds, as well as (in addition)
\beq \label{102}
\mbox{
$T$ is Hyers-Ulam $\Phi$-stable:\
$d(x,z)\le \Phi(d(x,Tx))$,\ for all $x\in X$.
}
\eeq
\etheorem

Note that, when $\vphi$ is {\it linear}
(i.e.: $\vphi(t)=\al t$, $t\in R_+$, for some $\al\in ]0,1[$)
the second conclusion above is necessarily retainable;
and Theorem \ref{t1} is just the 1922 
Banach's contraction mapping principle \cite{banach-1922}.
This result found some interesting applications in 
operator equations theory;
so, it was the subject of various extensions.
For example, a way of doing this
is by taking implicit "functional" contractive conditions 
\ben
\item[] (a02)\ \ 
$F(d(Tx,Ty),d(x,y))\le 0$,\ $\forall x,y\in X$, $x\ne y$;
\een
where $F:R_+^2\to R\cup\{-\oo,\oo\}$ is an appropriate function.
For more details about the possible choices of $F$
we refer to the 1976 paper by
Turinici \cite{turinici-1976}.
A further extension of (a02) is by considering
implicit functional contractive conditions like
\ben
\item[] (a03)\ \ 
$F(d(Tx,Ty),d(x,y),d(x,Tx),d(y,Ty),d(x,Ty),d(y,Tx))\le 0$,\\ 
$\forall x,y\in X$, $x\ne y$;
\een
where $F:R_+^6\to R\cup\{-\oo,\oo\}$ is an appropriate function.
Some concrete examples may be found in 
Berinde and Vetro \cite{berinde-vetro-2012};
see also
Rhoades \cite{rhoades-1977}
and
Akkouchi \cite{akkouchi-2010}.

Recently, an interesting contractive condition 
of the type (a02) was introduced in 
Wardowski \cite{wardowski-2012}.
It is our aim in the following to show that, 
concerning the fixed point question, 
a reduction to Theorem \ref{t1} is possible
for most of these contractions; cf. Section 3.
For the remaining ones, we provide in Section 4 
a result where some specific requirements posed by Wardowski
(cf. Section 5) are shown to be superfluous.
Finally, Section 2 has an introductory character.
Further aspects will be delineated in a future paper.

\section{Preliminaries}
\setcounter{equation}{0}

In the following, the concept of 
{\it semi-Wardowski} and {\it Wardowski} function are 
introduced; and some elementary facts 
about these are discussed. 

{\bf (A)}
Call $F:R_+\to R\cup\{-\oo\}$, 
a {\it semi-Wardowski} function, provided
\ben
\item[] (b01)\ \ 
$F(t)=-\oo$ if and only if $t=0$
\item[] (b02)\ \ 
$F$ is strictly increasing: $t< s$ $\limpl$ $F(t)< F(s)$.
\een
As a consequence of these facts, the lateral limits
\beq \label{201}
F(t-0):=\lim_{s\to t-}F(s),\
F(t+0):=\lim_{s\to t+}F(s)
\eeq
exist, for each $t> 0$; in addition,
\beq \label{202}
-\oo< F(t-0)\le F(t)\le F(t+0)< \oo,\  \forall t> 0.
\eeq
Note that, in general, $F$ is not continuous.
However, by the general properties of 
the  monotone functions, we have
(cf. Natanson, \cite[Ch 8, Sect 1]{natanson-1964}):

\bprop \label{p1}
Suppose that  $F:R_+\to R\cup\{-\oo\}$
is a semi-Wardowski function. Then,
there exists a countable subset 
$\De(F)\incl R_+^0:=]0,\oo[$, such that
\beq \label{203}
\mbox{
$F(t-0)=F(t)=F(t+0)$,\ for each $t\in R_+^0\sm \De(F)$.
}
\eeq
\eprop

In addition, the following property of such objects
is useful for us.

\blemma\label{le1}
Let $F:R_+\to R\cup\{-\oo\}$
be a semi-Wardowski function. Then,
\beq \label{204}
\forall t,s\in R_+:\ F(t)< F(s) \limpl t< s.
\eeq
\elemma

\bproof
Take the numbers $t,s\in R_+$ according to 
the premise of this relation. 
Clearly, $s> 0$; and, from this, the case $t=0$ 
is proved; hence, we may assume that $t> 0$.
The alternative $t=s$ gives $F(t)=F(s)$;
contradiction.
Likewise, the alternative $t> s$ yields
(via (b02)) $F(t)> F(s)$; again a contradiction.
Hence, $t< s$; and the conclusion follows.
\eproof

{\bf (B)}
Now, let us add one more condition
(upon such functions)
\ben
\item[] (b03)\ \ 
$F(t)\to -\oo$,\ when $t\to 0+$;\ i.e.: $F(0+)=-\oo$.
\een
When $F$ satisfies (b01)-(b03), it will be 
referred to as a {\it Wardowski} function.

\blemma \label{le2}
Suppose that  $F:R_+\to R\cup\{-\oo\}$
is a Wardowski function. Then,
for each sequence $(t_n)$ in $R_+^0$,
\beq \label{205}
\mbox{
$F(t_n)\to -\oo$ implies $t_n\to 0$.
}
\eeq
\elemma

\bproof
Suppose that this is not true:
there must be some $\veps> 0$ such that
\ben
\item[]
for each $n$, there exists $m\ge n$, such that :\ $t_m> \veps$.
\een
We get therefore a subsequence $(s_n:=t_{i(n)})$ 
of $(t_n)$ such that 
$$
\mbox{
$s_n> \veps$\ (hence, $F(s_n)> F(\veps)$),\ $\forall n$.
}
$$
This, however, contradicts the property $F(s_n)\to -\oo$;
hence the conclusion.
\eproof

\brem \label{r1}
\rm

From the above developments, it follows that 
the increasing property of $F$ may be 
taken as non-strict; i.e. (b02) may be replaced by
\ben
\item[] (b04)\ \ 
$F$ is increasing: $t\le s$ implies $F(t)\le F(s)$.
\een
However, we preferred to leave this "strict" version, 
for a better comparison with the original Wardowski's results.
\erem

\section{Left-continuous Wardowski functions}
\setcounter{equation}{0}

In the following, some auxiliary facts about 
left-continuous Wardowski functions are given.
Further, the first main 
result of this exposition is given.

{\bf (A)}
Suppose that  $F:R_+\to R\cup\{-\oo\}$ is a Wardowski function;
and let $a> 0$ be fixed in the sequel. 
Denote, for each $t\ge 0$,
\ben
\item[] (c01)\ \ 
$M(a,F)(t)=\{s\ge 0; a+F(s)\le F(t)\}$,\ 
$\vphi(t)=\sup M(a,F)(t)$.
\een
Note that, from (b03), $M(a,F)(t)$ is nonempty, 
for each $t\ge 0$; 
and $M(a,F)(0)=\{0\}$ [hence $\vphi(0)=0$].
In addition, for each $t> 0$ 
one has, via Lemma \ref{le1},
$$
s\in M(a,F)(t)\limpl F(s)< F(t)\limpl s< t; 
$$
wherefrom, again for all $t> 0$,
\beq \label{301}
\mbox{
$M(a,F)(t)\incl [0,t[$;\ whence, $\vphi(t)\le t$.
}
\eeq
Moreover, as $F$ is (strictly) increasing,
\beq \label{302}
0\le t_1\le t_2 \limpl M(a,F)(t_1)\incl M(a,F)(t_2)
\limpl  \vphi(t_1)\le \vphi(t_2);
\eeq
so that, $\vphi$ is increasing on $R_+$.

A basic problem to be posed is that of
determining sufficient conditions under $F$ such that
the (increasing) function $\vphi$ be regressive 
and Matkowski (tele-) admissible (cf. Section 1).
The following conditions will be 
taken into consideration here:

{\bf 3a)}
The Wardowski function $F:R_+\to R\cup\{-\oo\}$ fulfills
\ben
\item[] (c02)\ \ 
$F$ is left-continuous:\ $F(t-0)=F(t)$, $\forall t> 0$
\een

{\bf 3b)}
The Wardowski function $F:R_+\to R\cup\{-\oo\}$
fulfills (for some $k> 0$)
\ben
\item[] (c03)\ \ 
$F$ is $k$-regular:\ $t^kF(t)\to 0$, as $t\to 0+$.
\een
When $k\in ]0,1[$ is generic here, we say that 
$F$ is {\it regular}.

\bprop \label{p2}
Let the general conditions above be in use. 
Then,

{\bf i)}
If the Wardowski function
$F:R_+\to R\cup\{-\oo\}$ 
is left-continuous, then the associated 
function $\vphi$ is 
regressive and Matkowski admissible

{\bf ii)}
If the Wardowski function
$F:R_+\to R\cup\{-\oo\}$ 
is left-continuous and regular, 
then the associated 
function $\vphi$ is 
regressive and Matkowski tele-admissible.
\eprop

\bproof
There are three steps to be passed.

{\bf Step 1.}
Let $t> 0$ be arbitrary fixed; and put $u:=\vphi(t)$.
If $u=0$, we are done; so, without loss, we may
assume that $u> 0$.
By the very definition of this number, 
$$
\mbox{
$a+F(s)\le F(t)$,\ for all $s\in [0,u[$.
}
$$
So, passing to limit as $(s\to u-$), one gets 
(as $F$ is left-continuous)
\beq \label{303}
\mbox{
$a+F(u)\le F(t)$;\ i.e.: $a+F(\vphi(t))\le t$.
}
\eeq
This in turn yields $F(\vphi(t))< F(t)$; 
whence $\vphi(t)< t$.

{\bf Step 2.}
Fix some $t> 0$; and let the iterative
sequence $(t_n)$ be given as
[$t_0=t$, $t_{n+1}=\vphi(t_n)$; $n\ge 0$].
If $t_h=0$, for some index $h> 0$, we are done;
so, without loss, one may assume that 
[$t_n> 0$, $\forall n$].
By (\ref{303}), we have
$$
\mbox{
$a+F(t_{n+1})\le F(t_n)$,\ 
[hence, $a\le F(t_n)-F(t_{n+1})$],\ $\forall n$.
}
$$
Adding the first $n$ inequalities, one gets
$$
\mbox{
$na\le F(t_0)-F(t_n)$\ 
[hence: $F(t_n)\le F(t_0)-na$],\ $\forall n$;
}
$$
so that, passing to limit as $n\to \oo$, one derives
$F(t_n)\to -\oo$. 
This, along with Lemma \ref{le2}, gives $t_n=\vphi^n(t)\to 0$
as $n\to \oo$.

{\bf Step 3.}
Let $k\in ]0,1[$ be
given by the regularity of $F$.
From the above relation,
\beq \label{304}
nat_n^k\le [F(t_0)-F(t_n)]t_n^k,\ \forall n.
\eeq
By the convergence property of $(t_n; n\ge 0)$ and
the $k$-regularity of $F$,
the limit in the right hand side
is zero; so,
given $\be> 0$, there exists a rank $i=i(\be)$ 
such that 
$$
[F(t_0)-F(t_n)]t_n^k\le \be,\ \forall n\ge i.
$$
Combining with (\ref{304}) yields
(after transformations)
$$
\mbox{
$t_n\le (\be/an)^{1/k}$,\ for all $n\ge i$.
}
$$
This, along with the convergence of the series $\sum_{n\ge 1} n^{-1/k}$
tells us that so is the series 
$\sum_n t_n=\sum_n \vphi^n(t)$.
The proof is complete.
\eproof

{\bf (B)}
Let $(X,d)$ be a metric space;
and $T\in \calf(X)$ be a selfmap of $X$.
Given the real number $a> 0$
and the Wardowski function $F:R_+\to R\cup\{-\oo\}$,
let us say that $T$ is {\it $(a,F)$-contractive}, if
\ben
\item[] (c04)\ \ 
$a+F(d(Tx,Ty))\le F(d(x,y))$,\ $\forall x,y\in X$, $x\ne y$.
\een

The first main result of this exposition is

\btheorem \label{t2}
Suppose that $T$ is $(a,F)$-contractive, 
for some $a> 0$ and some Wardowski function $F$.
In addition, let $(X,d)$ be complete.
Then, 

{\bf j)}
If, in addition, $F$ is left-continuous, then
$T$ is a globally  strong Picard  operator (modulo $d$)

{\bf jj)}
If, moreover, 
$F$ is left-continuous and regular, then
$T$ is a globally strong tele-Picard operator (modulo $d$).
\etheorem

\bproof
By Proposition \ref{p2}, the associated 
function $\vphi\in \calf(in)(R_+)$ is
regressive and Matkowski admissible
(resp., Matkowski tele-admissible).
On the other hand, 
by (c04) and the very definition of $\vphi$, 
it results that 
$T$ is $(d;\vphi)$-contractive.
This, along with Theorem \ref{t1}, 
gives us all conclusions we need.
\eproof

\brem \label{r2}
\rm

All examples in 
Wardowski \cite{wardowski-2012}, like
\ben
\item[] (c05)\ \ 
$F(0)=-\oo$; $F(t)=\log(\al t^2+\be t)+\ga t$, $t> 0$,
\item[] (c06)\ \ 
$F(0)=-\oo$; $F(t)=(-1)(t^{-\de})$, $t> 0$,
\een
are illustrations of Theorem \ref{t2} above. 
[Here, $\log$  is the {\it natural logarithm} and
$\al,\be,\ga,\de > 0$ are constants].
Precisely, (c05) and the alternative $0< \de< 1$ 
of (c06)
may be handled with the second half of 
Theorem \ref{t2}.
On the other hand, the alternative $\de\ge 1$
of (c06)
is reducible to the first half of the same.
Finally, by a preceding observation, non-strict increasing
versions (modulo $F$) of all these are allowed;
we do not give further details.
\erem

\section{Discontinuous Wardowski functions}
\setcounter{equation}{0}

Let us now return to the general setting above.
As results from Proposition \ref{p2}, 
the left continuity requirement upon $F$ is essential
in Theorem \ref{t2} so as to deduce it 
(via Proposition \ref{p2}) 
from Theorem \ref{t1}.
In the absence of this, the reduction argument above
does not work.
And then, the question arises of 
to what extent is Theorem \ref{t2}
retainable (via its first half) 
in such an extended setting.
Strange enough, a positive answer to this is still available.
But, before stating it, some preliminaries are needed.

{\bf (A)}
Let $(X,d)$ be a metric space.
Call the sequence $(x_n; n\ge 0)$,
{\it $d$-semi-Cauchy}, when 
$d(x_n,x_{n+1})\to 0$ as $n\to \oo$.
Note that any $d$-Cauchy sequence is $d$-semi-Cauchy;
but the reciprocal is not in general true.
Concerning this last aspect, a useful property is available 
for such sequences which are not $d$-Cauchy.
Given the sequence $(r_n; n\ge 0)$ in $R$ and the point
$r\in R$, let us write 
\ben
\item[]
$r_n\to r+$ when [$r_n> r$, $\forall n$] and $r_n\to r$.
\een

\bprop \label{p3}
Suppose that $(x_n; n\ge 0)$ is $d$-semi-Cauchy but not $d$-Cauchy.
Further, let $\De$ be a countable part of $R_+^0$.
There exist then a number $\eta\in R_+^0\sm \De$, a rank $j(\eta)\ge 0$, 
and a couple
of rank-sequences $(m(j); j\ge 0)$, $(n(j); j\ge 0)$ with
\beq \label{401}
j\le m(j)< n(j),\ d(x_{m(j)},x_{n(j)})> \eta,\ 
\forall j\ge 0
\eeq
\beq \label{402}
n(j)-m(j)\ge 2,\ d(x_{m(j)},x_{n(j)-1})\le \eta,\ 
\forall j\ge j(\eta)
\eeq
\beq \label{403}
d(x_{m(j)},x_{n(j)})\to \eta+,\ \mbox{as}\ j\to \oo
\eeq
\beq \label{404}
d(x_{m(j)+p},x_{n(j)+q})\to \eta,\ \mbox{as}\ j\to \oo,\ 
\forall p,q\in \{0,1\}.
\eeq
\eprop

\bproof
As $R_+^0\sm \De$ is dense in $R_+^0$, 
the $d$-Cauchy property of $(x_n; n\ge 0)$ writes
\ben
\item[]
$\forall \eta\in R_+^0\sm \De$, $\exists k=k(\eta)$:\
$k\le m< n$ $\limpl$ $d(x_m,x_n)\le \eta$.
\een
The negation of this property means:
there must be some $\eta\in R_+^0\sm \De$,
such that
\beq \label{405}
(\forall j\ge 0):\ 
A(j):=\{(m,n)\in N\times N; j\le m< n, d(x_m,x_n)> \eta\}
\ne \es.
\eeq
Having this precise, denote, for each $j\ge 0$,
\ben
\item[] (d01)\ \ 
$m(j)=\min \Dom(A(j))$,\ $n(j)=\min A(m(j))$.
\een
As a consequence, the couple of rank-sequences 
$(m(j); j\ge 0)$, $(n(j); j\ge 0)$  fulfills (\ref{401}).
On the other hand, letting
the index $j(\eta)$ be such that
\beq \label{406}
\mbox{
$\rho_i:=d(x_i,x_{i+1})< \eta$,\ for all $i\ge j(\eta)$,
}
\eeq
it is clear that  (\ref{402}) holds too.
This in turn yields (by the triangular inequality)
$$  \barr{l}
\eta< d(x_{m(j)},x_{n(j)})\le 
d(x_{m(j)},x_{n(j)-1})+\rho_{n(j)-1}\le 
\eta+\rho_{n(j)-1},\ \forall j\ge j(\eta);
\earr
$$
so, passing to limit as $j\to \oo$ gives (\ref{403}).
Finally, again by the triangular inequality,
$$
(\forall j\ge 0):\
d(x_{m(j)},x_{n(j)})-\rho_{n(j)}\le
d(x_{m(j)},x_{n(j)+1})\le 
d(x_{m(j)},x_{n(j)})+\rho_{n(j)}.
$$
By a limit process upon $j$,  one gets 
the case $(p=0,q=1)$ of (\ref{404}).
The remaining ones are obtained in a similar way.
\eproof

{\bf (B)}
We are now in position to state
the second main result of this exposition.
Let $T\in \calf(X)$ be a selfmap of $X$.

\btheorem \label{t3}
Suppose that $T$ is $(a,F)$-contractive, 
for some $a> 0$ and some Wardowski function 
$F:R_+\to R\cup\{-\oo\}$.
In addition, let $(X,d)$ be complete.
Then, $T$ is a globally  strong Picard 
operator (modulo $d$).
\etheorem

\bproof
From the $(a,F)$-contractive condition,
it results that
\beq \label{407}
\mbox{
$T$ is a strict contraction:\ 
$d(Tx,Ty)< d(x,y)$,\ $\forall x,y\in X$, $x\ne y$.
}
\eeq
This firstly gives us that $\Fix(T)$ is an asingleton. 
As a second consequence, 
\beq \label{408}
\mbox{
$T$ is nonexpansive:\ $d(Tx,Ty)\le d(x,y)$,\ $\forall x,y\in X$;
}
\eeq
hence, in particular, $T$ is $d$-continuous on $X$.
So, it remains only to prove that $T$ is Picard 
(modulo $d$).
Given $x=x_0$ in $X$, 
put ($x_n:=T^nx_0$; $n\ge 0$).
If $x_i=x_{i+1}$ for some $i\ge 0$, we are done;
so, without loss, one may assume that
\ben
\item[] (d02)\ \ 
$x_n\ne x_{n+1}$ (hence, $\rho_n:=d(x_n,x_{n+1})> 0$),\ $\forall n$.
\een

{\bf Part 1}.
By the contractive condition, we have
$$
\mbox{
$a+F(\rho_{n+1})\le F(\rho_n)$,\ 
[hence, $a\le F(\rho_n)-F(\rho_{n+1})$],\ $\forall n$.
}
$$
Adding the first $n$ inequalities, one gets
$$
\mbox{
$na\le F(\rho_0)-F(\rho_n)$\ 
[hence: $F(\rho_n)\le F(\rho_0)-na$],\ $\forall n$;
}
$$
so that, passing to limit as $n\to \oo$, one derives
$F(\rho_n)\to -\oo$. 
This, along with Lemma \ref{le2}, gives $\rho_n\to 0$;
hence, $(x_n; n\ge 0)$ is $d$-semi-Cauchy.

{\bf Part 2}.
Let $\De(F)$ stand for the subset of
$R_+^0$ where $F$ is discontinuous;
note that, by Proposition \ref{p1}, 
it is (at most) countable.
Assume by contradiction that $(x_n)$ is not $d$-Cauchy.
By Proposition \ref{p3}, there exists
a number $\eta\in R_+^0\sm \De(F)$, a rank $j(\eta)\ge 0$, 
and a couple
of rank-sequences $(m(j); j\ge 0)$, $(n(j); j\ge 0)$ with
the properties (\ref{401})-(\ref{404}).
By the contractive condition, we have
$$
a+F(d(x_{m(j)+1},x_{n(j)+1}))\le F(d(x_{m(j)},x_{n(j)})),\ \forall j.
$$
Passing to limit as $j\to \oo$ one gets,
from the choice of $\eta$ (as a point where $F$ is 
(bilaterally) continuous) and (\ref{403})-(\ref{404}),
$$
\mbox{
$a+F(\eta)\le F(\eta)$;\ hence, $a\le 0$.
}
$$
The obtained contradiction tells us that 
$(x_n)$ is $d$-Cauchy; and, from this,
all is clear via (\ref{408}).
The proof is complete.
\eproof

\section{Further aspects}
\setcounter{equation}{0}

Let again $(X,d)$ be a metric space;
and $T\in \calf(X)$ be a selfmap of $X$;
supposed to be $(a,F)$-contractive, 
for some $a> 0$ and some Wardowski function 
$F:R_+\to R\cup\{-\oo\}$.
According to Theorem \ref{t3},
$T$ is then a globally strong Picard operator 
(modulo $d$).
Under such a perspective,
any supplementary conditions upon these data
is superfluous if we want that 
{\it exactly} this property be reached.
However, for the property described by
the second half of Theorem \ref{t2} being available,
the question of identifying these
conditions is not at all superfluous.
As suggested by the quoted result, 
an extra condition to be considered is 
the regularity of $F$.
This is, indeed, in effect for our purpose;
as certified by

\btheorem \label{t4}
Suppose that $T$ is $(a,F)$-contractive, 
for some $a> 0$ and some regular 
Wardowski function $F$.
In addition, let $(X,d)$ be complete.
Then, $T$ is a globally  strong tele-Picard 
operator (modulo $d$).
\etheorem

This result is, essentially, the one in
Wardowski \cite{wardowski-2012}.
For completeness reasons, we shall 
provide its proof, with some modifications.

\bproof {\bf (Theorem \ref{t4})}
As in Theorem \ref{t3}, we only have to 
establish that $T$ is tele-Picard (modulo $d$).
Given $x=x_0$ in $X$, put $(x_n:=T^nx_0; n\ge 0)$.
If $x_i=x_{i+1}$ for some $i\ge 0$, we are done;
so, without loss, one may assume that 
[$x_n\ne x_{n+1}$ (hence, $\rho_n:=d(x_n,x_{n+1})> 0$), $\forall n$].
By the contractive condition, we have
$$
\mbox{
$a+F(\rho_{n+1})\le F(\rho_n)$,\ 
[hence, $a\le F(\rho_n)-F(\rho_{n+1})$],\ $\forall n$.
}
$$
Adding the first $n$ inequalities, one gets
\beq \label{501}
\mbox{
$na\le F(\rho_0)-F(\rho_n)$\ 
[hence: $F(\rho_n)\le F(\rho_0)-na$],\ $\forall n$;
}
\eeq
This gives (passing to limit as $n\to \oo$),
$F(\rho_n)\to -\oo$;  
wherefrom, by Lemma \ref{le2}, $\rho_n\to 0$;
hence, $(x_n; n\ge 0)$ is $d$-semi-Cauchy.
Further, let $k\in ]0,1[$ be given by the regularity of $F$.
Again via (\ref{501}), one gets
\beq \label{502}
na\rho_n^k\le [F(\rho_0)-F(\rho_n)]\rho_n^k,\ \forall n.
\eeq
By the convergence property of $(\rho_n; n\ge 0)$ and
the $k$-regularity of $F$, the limit of the right hand side
is zero; so that,
given $\be> 0$, there exists $i=i(\be)$ 
such that 
$$
[F(\rho_0)-F(\rho_n)]\rho_n^k\le \be,\ \forall n\ge i.
$$
Combining with (\ref{502}) yields
(after transformations)
$$
\mbox{
$\rho_n\le (\be/an)^{1/k}$,\ for all $n\ge i$.
}
$$
This, along with the convergence of the series $\sum_{n\ge 1} n^{-1/k}$
tells us that so is the series 
$\sum_n \rho_n=\sum_n d(x_n,x_{n+1})$; wherefrom, $(x_n; n\ge 0)$ is $d$-Cauchy.
The last part is identical with the one of
Theorem \ref{t3}; and conclusion follows.
\eproof

Note, finally,  that all these results are extendable
to the framework of quasi-ordered metric spaces
under the lines in 
Agarwal et al \cite{agarwal-el-gebeily-o-regan-2008};
see also
Turinici \cite{turinici-1986}.
A development of these facts will be given in a future paper.


\end{document}